# A Sliding Blocks Estimator for the Extremal Index

**Christian Y. Robert**

*Ecole Nationale de la Stastistique et de l'Administration Economique, Timbre J120*
*3 Avenue Pierre Larousse*
*92245 Malakoff Cedex, France*
*e-mail:* chrobert@ensae.fr

**Johan Segers***

*Institut de statistique, Université catholique de Louvain*
*Voie du Roman Pays 20*
*B-1348 Louvain-la-Neuve, Belgium*
*e-mail:* johan.segers@uclouvain.be

**Christopher A. T. Ferro**

*School of Engineering, Computing and Mathematics, University of Exeter*
*Harrison Building, North Park Road*
*Exeter, EX4 4QF, United Kingdom*
*e-mail:* c.a.t.ferro@exeter.ac.uk

**Abstract:** In extreme value statistics for stationary sequences, blocks estimators are usually constructed by using disjoint blocks because exceedances over high thresholds of different blocks can be assumed asymptotically independent. In this paper we focus on the estimation of the extremal index which measures the degree of clustering of extremes. We consider disjoint and sliding blocks estimators and compare their asymptotic properties. In particular we show that the sliding blocks estimator is more efficient than the disjoint version and has a smaller asymptotic bias. Moreover we propose a method to reduce its bias when considering sufficiently large block sizes.

**AMS 2000 subject classifications:** Primary 60G70, 62E20; secondary 62G20, 62G32.
**Keywords and phrases:** clusters of extremes, extremal index, FTSE 100, intervals estimator, max-autoregressive process, moving maximum process, maximal correlation coefficient, mixing coefficient, sample maximum, stationary time series.

## Contents

*Research supported by IAP research network grant nr. P6/03 of the Belgian government (Belgian Science Policy).

## 1. Introduction

Suppose that $(X_n)_{n\in\mathbb{N}}$ is a strictly stationary sequence of random variables with marginal distribution function $F$. We assume that this sequence has an extremal index $\theta \in (0,1]$, that is, for each $\tau > 0$, there exists a sequence of levels $(u_n(\tau))_{n\in\mathbb{N}}$ such that $\lim_{n\to\infty} n\bar{F}\big(u_n(\tau)\big) = \tau$ and

$$\lim_{n\to\infty} \Pr\big(M_n \leqslant u_n(\tau)\big) = e^{-\theta\tau},$$

where $\bar{F} = 1 - F$ and $M_n = \max\{X_1, \ldots, X_n\}$. The extremal index can be interpreted in a number of ways, the most common one being the reciprocal of the mean cluster size in the limiting point process of exceedance times over high thresholds. The probabilistic theory was worked out in (11), (13), (15), (9), (14), and (12).

Our objective is to estimate $\theta$ based on a finite stretch $X_1, \ldots, X_n$ from the time series. Inference about the extremal index parameter has been extensively studied. The three more common approaches are the blocks method, the runs method and the inter-exceedance times method. The two first methods identify clusters and construct estimates for $\theta$ based on these clusters. For each of these methods, there are two parameters which determine the clusters and consequently the estimates of $\theta$: a threshold and a cluster identification scheme parameter. The third method is based on inter-exceedance times and obviates the need for a cluster identification scheme parameter. Some references on estimation of the extremal index using these three approaches are (7), (8), (17), (18), (6), (10) and (16) among others.

In this paper we focus on the blocks method. Traditionally, it consists of partitioning the $n$ observations into consecutive blocks of a certain length, say $r$. In each block, the number of exceedances over a certain high threshold are counted, and the blocks estimator is then defined as the reciprocal of the average number of exceedances per block among blocks with at least one exceedance.

Blocks estimators are usually constructed by using disjoint blocks, for in that case the blocks can be assumed to be approximately independent.

The main novelty in this paper is our proposal to use sliding rather than disjoint blocks, that is, to slide a window of length $r$ through the sample, yielding $n - r + 1$ blocks rather than just $n/r$ disjoint blocks. Surprisingly, this simple modification leads to a more efficient estimator with a smaller asymptotic variance. Moreover we provide estimators of the asymptotic variances of the estimators, which permits the construction of confidence intervals and the selection of variance-minimizing thresholds. We also provide a way to estimate and correct for the asymptotic bias of the estimators.

In contrast to most previous papers but in accordance with (16), we assume that thresholds and block sizes are such that the expected number of excesses per block converges to a positive constant. In practice, the threshold is chosen as a large order statistic. However, mathematical treatment of such random thresholds requires complicated empirical process techniques.

The content of the paper is organized as follows. In Section 2 we introduce the blocks estimators for the extremal index. In Section 3 we consider asymptotic variances and covariances of the mean number of excesses per block and the empirical distribution functions of disjoint and sliding block maxima. We establish consistency and asymptotic normality of our estimators in Section 4. We discuss how to estimate and minimize their asymptotic variance in Section 5 and how to reduce their bias in Section 6. In Section 7, we investigate the finite sample behavior of the estimators on simulated data and we provide a case study. Proofs are spelled out in the appendices.

## 2. The estimators

For positive integer $r$, put

$$F_r(u) := \Pr(M_r \leqslant u), \qquad \tau_r(u) := r\bar{F}(u) \qquad \text{and} \qquad \theta_r(u) := -\frac{\log F_r(u)}{\tau_r(u)}. \tag{2.1}$$

It follows from the definition of the extremal index that $\theta = \lim_{r\to\infty} \theta_r(u_r)$ where $u_r = u_r(\tau)$. The estimators of $\theta$ to be proposed are based upon empirical analogues of the functions $F_r$ and $\tau_r$.

For integer $0 \leqslant s < r$, put $M_{s,r} := \max_{s<i\leqslant r} X_i$. Note that $M_r = M_{0,r}$. The distribution function $F_r$ of the block maximum $M_r$ can be estimated using maxima of $k := \lfloor n/r \rfloor$ *disjoint* blocks or using maxima of $n-r+1$ *sliding* blocks:

$$\hat{F}^{\mathrm{dj}}_{n,r}(u) := \frac{1}{k}\sum_{i=1}^{k} I(M_{(i-1)r,ir} \leqslant u), \quad \hat{F}^{\mathrm{sl}}_{n,r}(u) := \frac{1}{n-r+1}\sum_{i=0}^{n-r} I(M_{i,i+r} \leqslant u).$$

One may wonder why the use of sliding rather than disjoint blocks should make a difference. After all, the $n - r + 1$ blocks in the definition of $\hat{F}^{\mathrm{sl}}_{n,r}(u)$ are overlapping and hence strongly dependent, even in the iid case. Nevertheless, we will show in Proposition 3.1 below that the asymptotic variance of $\hat{F}^{\mathrm{sl}}_{n,r}(u)$ is typically smaller than the one of $\hat{F}^{\mathrm{dj}}_{n,r}(u)$.

Writing

$$\hat{\tau}_{n,r}(u) := \frac{1}{k}\sum_{j=1}^{k}\sum_{i=1}^{r} I(X_{r(j-1)+i} > u) = \frac{1}{k}\sum_{i=1}^{rk} I(X_i > u) \tag{2.2}$$

the *disjoint* and the *sliding blocks estimators* of the extremal index can now be defined as follows:

$$\hat{\theta}^{\mathrm{dj}}_{n,r}(u) := -\frac{\log \hat{F}^{\mathrm{dj}}_{n,r}(u)}{\hat{\tau}_{n,r}(u)}, \qquad \hat{\theta}^{\mathrm{sl}}_{n,r}(u) := -\frac{\log \hat{F}^{\mathrm{sl}}_{n,r}(u)}{\hat{\tau}_{n,r}(u)}.$$

As above, the sliding version will turn out to be more efficient than the disjoint one, see Corollary 4.3.

The estimators require the choice of two tuning parameters: the threshold $u$ and the block size $r$. If $u$ is equal to the $\lfloor k\tau \rfloor$-th largest order statistic of $X_1, \ldots, X_n$, the disjoint blocks estimator is the same as $\hat{\theta}^{(\tau)}_{n,1}$ in (16). As mentioned in Section 1, the mathematical treatment of such random thresholds is intricate and requires empirical process techniques. For the sake of simplicity, the threshold sequence $(u_r)_{r\in\mathbb{N}}$ will be assumed to be deterministic. Comparing our Corollary 4.3 with **(author?)** (16, Corollary 4.2), it follows that this simplifying assumption does not make any difference asymptotically.

## 3. Asymptotic variances and covariances

The disjoint and sliding blocks estimators for $\theta$ are functions of $\hat{F}^{\mathrm{dj}}_{n,r}(u)$, $\hat{F}^{\mathrm{sl}}_{n,r}(u)$, and $\hat{\tau}_{n,r}(u)$. We shall need to find the asymptotic variances and covariances of the latter three estimators. Most importantly, we will show that $\hat{F}^{\mathrm{dj}}_{n,r}(u) - \hat{F}^{\mathrm{sl}}_{n,r}(u)$ has a non-negligible asymptotic variance and is asymptotically uncorrelated with $\hat{F}^{\mathrm{sl}}_{n,r}(u)$. As a result, the sliding blocks estimator for $F_r(u)$ is the most efficient convex combination of the disjoint and sliding blocks estimators for $F_r(u)$. The proofs of the results in this section are to be found in Appendix A.

The *maximal correlation coefficients* of the process $(X_n)_{n\in\mathbb{N}}$ relative to the threshold $u$ are defined by

$$\rho_{n,l}(u) := \max_{t=1,\ldots,n-l} \sup_{\substack{\xi \in L_2(\mathcal{F}_{1,t}(u)) \\ \eta \in L_2(\mathcal{F}_{t+l,n}(u))}} |\operatorname{corr}(\xi,\eta)|.$$

Here $\mathcal{F}_{a,b}(u)$ is the $\sigma$-field generated by the events $\{X_i \leqslant u\}$ for $i \in \{a, \ldots, b\}$, and $L_2(\mathcal{F})$ is the space of $\mathcal{F}$-measurable square-integrable random variables. Obviously, the random variables $\xi$ and $\eta$ in the definition of $\rho_{n,l}(u)$ should have positive variance. For comparisons with other mixing coefficients, see e.g. (2). Here we just wish to note that

$$\rho_{n,l}(u) \geqslant \alpha_{n,l}(u) := \max_{t=1,\ldots,n-l} \sup_{\substack{A \in \mathcal{F}_{1,t}(u) \\ B \in \mathcal{F}_{t+l,n}(u)}} |\Pr(A\cap B) - \Pr(A)\Pr(B)|.$$

The coefficients $\alpha_{n,l}$ underlie the condition called $\Delta(u_n)$ in (9) and are themselves greater than the coefficients introduced in (11) yielding Leadbetter's $D(u_n)$ condition. Since the upper bounds we will impose on $\rho_{n,l}$ will trivially imply the same upper bounds on $\alpha_{n,l}$, the results in (9) become available to us as well.

Let $r_n$ and $l_n$ be positive integer sequences such that, as $n \to \infty$,

$$l_n = o(r_n), \qquad r_n = o(n) \qquad \text{and} \qquad \sum_{l=l_n}^{n} \rho_{n,l}(u_{r_n}) = o(r_n). \tag{3.1}$$

Note that the assumptions imply $\sum_{l=1}^{n} \rho_{n,l}(u_{r_n}) = o(r_n)$ and that the final assumption is implied by $k_n \rho_{n,l_n}(u_{r_n}) \to 0$ as $n \to \infty$, where $k_n = \lfloor n/r_n \rfloor$.

**Proposition 3.1.** *Let $(X_n)_{n\in\mathbb{N}}$ be stationary with extremal index $\theta$ and let $(u_r)_{r\in\mathbb{N}}$ be a sequence of thresholds such that $r\bar{F}(u_r) \to \tau \in (0,\infty)$ as $r \to \infty$. If (3.1) holds, then, as $n \to \infty$, denoting $\alpha := \theta\tau$,*

$$\begin{aligned} k_n \operatorname{var}\big(\hat{F}^{\mathrm{dj}}_{n,r_n}(u_{r_n})\big) &\to e^{-\alpha}(1-e^{-\alpha}), &(3.2)\\ k_n \operatorname{var}\big(\hat{F}^{\mathrm{sl}}_{n,r_n}(u_{r_n})\big) &\to 2\alpha^{-1}e^{-\alpha}\big(1-(1+\alpha)e^{-\alpha}\big), &(3.3)\\ k_n \operatorname{cov}\big(\hat{F}^{\mathrm{sl}}_{n,r_n}(u_{r_n}), \hat{F}^{\mathrm{dj}}_{n,r_n}(u_{r_n})\big) &\to 2\alpha^{-1}e^{-\alpha}\big(1-(1+\alpha)e^{-\alpha}\big). &(3.4) \end{aligned}$$

By Proposition 3.1, we can write $\hat{F}^{\mathrm{dj}}_{n,r_n}(u_{r_n}) = \hat{F}^{\mathrm{sl}}_{n,r_n}(u_{r_n}) + \varepsilon_n$, the random term $\varepsilon_n$ having mean zero, being asymptotically uncorrelated with $\hat{F}^{\mathrm{sl}}_{n,r_n}(u_{r_n})$, and having non-negligible asymptotic variance. As a consequence, the sliding blocks estimator of the distribution function of the block maximum is more efficient than the disjoint blocks estimator. The two asymptotic variance functions as well as their ratio are shown in Figure 1. Observe that the relative efficiency of the disjoint versus the sliding blocks estimator is decreasing in $\alpha$. For $\alpha \to 0$, the clusters of exceedances become very sparse, and the two estimators are asymptotically equivalent.

In order to get the asymptotic covariances between $\hat{\tau}_{n,r_n}(u_{r_n})$ in (2.2) with the disjoint and sliding blocks estimators for $F_r(u)$, a somewhat stronger condition on the maximal correlation coefficients is needed:

$$\sum_{l=1}^{n} \rho_{n,l}(u_{r_n}) = o(r_n^{1/2}), \qquad n \to \infty. \tag{3.5}$$

**Proposition 3.2.** *If in addition to the conditions of Proposition 3.1 also (3.5) holds, then as $n \to \infty$,*

$$\left.\begin{aligned} &k_n \operatorname{cov}\big(\hat{F}^{\mathrm{dj}}_{n,r_n}(u_{r_n}), \hat{\tau}_{n,r_n}(u_{r_n})\big)\\ &k_n \operatorname{cov}\big(\hat{F}^{\mathrm{sl}}_{n,r_n}(u_{r_n}), \hat{\tau}_{n,r_n}(u_{r_n})\big) \end{aligned}\right\} \to -\tau e^{-\alpha}. \tag{3.6}$$

Finally, in order to find the asymptotic variance of $\hat{\tau}_{n,r_n}(u_{r_n})$, another additional assumption is needed: there exists a positive integer sequence $(s_n)_{j\in\mathbb{N}}$ and a probability distribution $(\pi_j)_{j\in\mathbb{N}}$ on the positive integers such that as $n \to \infty$

$$l_n = o(s_n), \qquad s_n = o(r_n), \qquad \rho_{n,l_n}(u_{r_n}) = o(s_n/n) \tag{3.7}$$

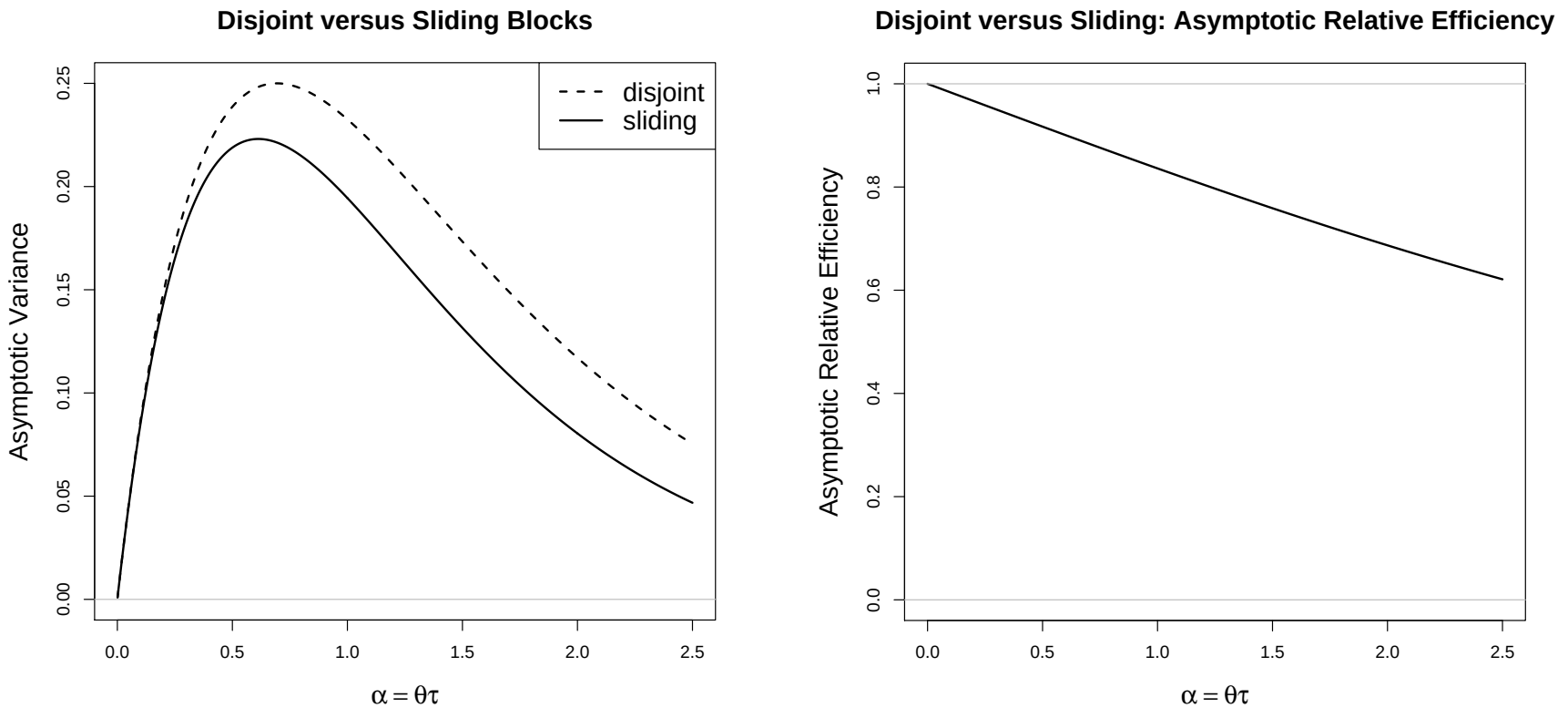

FIG 1. *Left:* Asymptotic variance functions of disjoint and sliding blocks estimators of $F_r(u)$, see Proposition 3.1. *Right:* Asymptotic relative efficiency of disjoint versus sliding blocks estimators of $F_r(u)$.

as well as, writing $N_s(u) = \sum_{i=1}^{s} I(X_i > u)$,

$$\begin{aligned} &\Pr\{N_{s_n}(u_{r_n}) = j \mid M_{s_n} > u_{r_n}\} \to \pi_j \quad \text{for all } j \in \mathbb{N}, \\ &\mathrm{E}[N_{s_n}^2(u_{r_n}) \mid M_{s_n} > u_{r_n}] \to \textstyle\sum_{j=1}^{\infty} j^2 \pi_j < \infty. \end{aligned} \tag{3.8}$$

The distribution $(\pi_j)_{j\in\mathbb{N}}$ is called the *cluster size distribution*; it describes the limiting probability distribution of the number of threshold excesses within the block $X_1, \ldots, X_{s_n}$ given that there is at least one such excess. The second part of (3.8) is a uniform integrability condition ensuring that the first two moments of the finite-sample cluster size distribution converge to the proper limits. Note that $\Pr(M_{s_n} > u_{r_n}) \leqslant s_n \Pr(X_1 > u_{r_n}) \to 0$ as $n \to \infty$ while

$$\frac{s_n \Pr(X_1 > u_{r_n})}{\Pr(M_{s_n} > u_{r_n})} = \mathrm{E}[N_{s_n}(u_{r_n}) \mid M_{s_n} > u_{r_n}] \to \sum_{j=1}^{\infty} j\pi_j = \theta^{-1}.$$

Under the above conditions, the asymptotic distribution of $N_{r_n}(u_{r_n})$ is compound Poisson (9, Theorem 5.1):

$$N_{r_n}(u_{r_n}) \xrightarrow{d} N := \sum_{i=1}^{\nu} \zeta_i \tag{3.9}$$

where $\nu$ is a Poisson$(\theta\tau)$ random variable and $(\zeta_i)_{i\in\mathbb{N}}$ is a sequence of positive independent and identically distributed integer-valued random variables from the cluster size distribution, independent of $\nu$. Note that $\mathrm{E}(\zeta_1) = \sum_{j\geqslant 1} j\pi_j = \theta^{-1}$. Moreover,

$$\mathrm{E}(N) = \tau, \qquad \mathrm{var}(N) = \alpha \sum_{j=1}^{\infty} j^2 \pi_j. \tag{3.10}$$

**Proposition 3.3.** *If in addition to the conditions of Proposition 3.1 also (3.7)–(3.8) hold, then as* $n \to \infty$,

$$k_n \operatorname{var}\bigl(\hat{\tau}_{n,r_n}(u_{r_n})\bigr) \to \operatorname{var}(N) = \alpha \sum_{j=1}^{\infty} j^2 \pi_j. \tag{3.11}$$

## 4. Weak consistency and asymptotic normality

The main result of this paper is the joint asymptotic normality of the disjoint and sliding blocks estimators for $\theta$ in Corollary 4.3. The proofs of the results in this section are to be found in Appendix B. Write

$$m_p := \sum_{j \geqslant 1} j^p \pi_j, \qquad p \in \{1, 2\}.$$

Recall that $m_1 = \theta^{-1}$.

**Theorem 4.1.** *Let* $(X_n)_{n\in\mathbb{N}}$ *be stationary with extremal index* $\theta$ *and let* $(u_r)_{r\in\mathbb{N}}$ *be a sequence of thresholds such that* $r\bar{F}(u_r) \to \tau \in (0,\infty)$ *as* $r \to \infty$. *If conditions (3.1), (3.7) and (3.8) hold, then as* $n \to \infty$

$$\hat{\theta}^{\mathrm{dj}}_{n,r_n}(u_{r_n}) \xrightarrow{p} \theta \qquad \text{and} \qquad \hat{\theta}^{\mathrm{sl}}_{n,r_n}(u_{r_n}) \xrightarrow{p} \theta.$$

In order to get asymptotic normality of the estimators, we will need an additional technical assumption: there exists a constant $p$ with $p > 1$ such that as $n \to \infty$,

$$\mathrm{E}[N^{2p}_{r_n}(u_{r_n})] = O(1). \tag{4.1}$$

We first state joint asymptotic normality of $\hat{F}^{\mathrm{dj}}_{n,r}(u)$, $\hat{F}^{\mathrm{sl}}_{n,r}(u)$, and $\hat{\tau}_{n,r}(u)$. Joint asymptotic normality of $\hat{\theta}^{\mathrm{dj}}_{n,r}(u)$ and $\hat{\theta}^{\mathrm{sl}}_{n,r}(u)$ then follows by the delta-method.

**Theorem 4.2.** *Let* $(X_n)_{n\in\mathbb{N}}$ *be stationary with extremal index* $\theta$ *and let* $(u_r)_{r\in\mathbb{N}}$ *be a sequence of thresholds such that* $r\bar{F}(u_r) \to \tau \in (0,\infty)$ *as* $r \to \infty$. *If conditions (3.1), (3.5), (3.7)–(3.8) and (4.1) hold, then*

$$\sqrt{k_n} \begin{pmatrix} \hat{F}^{\mathrm{dj}}_{n,r_n}(u_{r_n}) - F_{r_n}(u_{r_n}) \\ \hat{F}^{\mathrm{sl}}_{n,r_n}(u_{r_n}) - F_{r_n}(u_{r_n}) \\ \hat{\tau}_{n,r_n}(u_{r_n}) - \tau_{r_n}(u_{r_n}) \end{pmatrix} \xrightarrow{d} N(\mathbf{0}, \mathbf{\Sigma})$$

*where* $\mathbf{\Sigma} = (\sigma_{ij})_{i,j=1}^3$ *is symmetric and, writing* $\alpha = \tau\theta$,

$$\begin{aligned} \sigma_{11} &= e^{-\alpha}(1 - e^{-\alpha}), \\ \sigma_{22} = \sigma_{12} &= 2\alpha^{-1} e^{-\alpha}\bigl(1 - (1+\alpha)e^{-\alpha}\bigr), \\ \sigma_{31} = \sigma_{32} &= -\tau e^{-\alpha}, \\ \sigma_{33} &= \alpha m_2. \end{aligned}$$

Recall $\theta_r(u)$ in (2.1). In order to control the bias of the extremal index estimators, assume that the block sizes are sufficiently large so that as $n \to \infty$,

$$\theta_{r_n}(u_{r_n}) - \theta = o\big(1/\sqrt{k_n}\big). \tag{4.2}$$

The asymptotic variance of the extremal index estimators will depend on $\theta$, $\tau$, and the squared coefficient of variation $c^2$ of the cluster size distribution $(\pi_j)_{j\in\mathbb{N}}$:

$$c^2 := \frac{m_2 - m_1^2}{m_1^2}.$$

**Corollary 4.3.** *Under the conditions of Theorem 4.2 and (4.2), as $n \to \infty$,*

$$\sqrt{k_n} \begin{pmatrix} \hat{\theta}^{\mathrm{dj}}_{n,r_n}(u_{r_n}) - \theta \\ \hat{\theta}^{\mathrm{sl}}_{n,r_n}(u_{r_n}) - \theta \end{pmatrix} \xrightarrow{d} N(\mathbf{0}, \boldsymbol{V})$$

*where $\boldsymbol{V} = (v_{ij})_{i,j=1}^2$ is symmetric and*

$$v_{11} = \frac{\theta^2}{\alpha^2}(e^\alpha - 1 - \alpha) + \frac{\theta^2 c^2}{\alpha},$$
$$v_{22} = v_{12} = 2\frac{\theta^2}{\alpha^3}\left(e^\alpha - 1 - \alpha - \frac{\alpha^2}{2}\right) + \frac{\theta^2 c^2}{\alpha}.$$

The asymptotic variance of the disjoint blocks estimator corresponds with the one for the same estimator but at a random threshold (order statistic) in **(author?)** (16, Corollary 4.2). It is worth noting that $v_{22} \leqslant v_{11}$. As a result, the sliding blocks estimator is more efficient than its disjoint version. Even more, the most efficient convex combination of the disjoint and sliding blocks estimators is the sliding blocks estimator itself.

## 5. Estimating and minimizing the asymptotic variance

For a fixed $c^2 \geqslant 0$, the asymptotic variance functions of $\sqrt{k_n}(\hat{\theta}_{n,r_n}/\theta - 1)$,

$$0 \leqslant \alpha \mapsto \begin{cases} \alpha^{-2}(e^\alpha - 1 - \alpha) + \alpha^{-1}c^2, & \text{disjoint blocks,} \\ 2\alpha^{-3}\left(e^\alpha - 1 - \alpha - \frac{\alpha^2}{2}\right) + \alpha^{-1}c^2, & \text{sliding blocks,} \end{cases} \tag{5.1}$$

are convex and possess unique global minima. These minima and the values of $\alpha$ for which they are attained can be computed numerically, see Figure 2. Hence, given an estimate of $c^2$, we can estimate the respective optimal values for $\alpha$, divide by an estimate of $\theta$, and thus obtain estimates of the optimal $\tau$ to be used for the disjoint or sliding blocks estimators. Given such estimates, we can for a given threshold $u$ estimate the asymptotically optimal block lengths $r$ and vice versa.

The missing element in this procedure is an estimate of $c^2$. Knowledge of $c^2$ is also needed when one wants to construct asymptotic confidence intervals for $\theta$ based on Corollary 4.3 or estimate the asymptotic bias of the extremal index estimators in Section 6 below. In addition, the quantity $c^2$ is interesting in its own right as a measure of dispersion of the cluster size distribution $(\pi_j)_{j\in\mathbb{N}}$. Since the mean cluster size is equal to $m_1 = \theta^{-1}$, for which consistent estimators are available, we can focus here on estimating the cluster-size variance $m_2 - m_1^2$ or the second moment $m_2$.

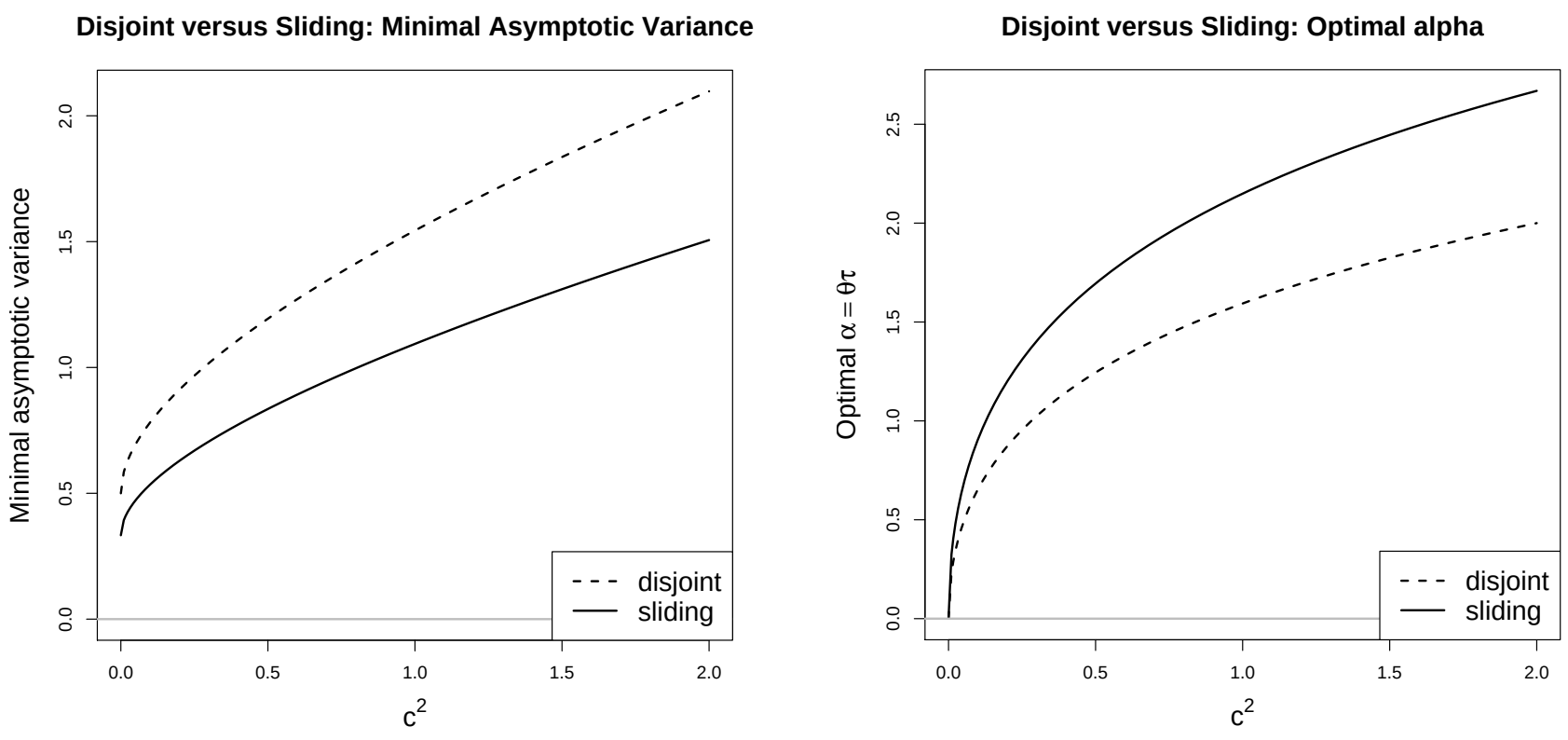

Fig 2. Minima (left) and minimizers (right) of the asymptotic variance functions in (5.1) of the disjoint and sliding blocks estimators for $\theta$.

A first possible strategy to estimate the cluster-size variance is to partition the threshold exceedances into clusters and estimate the cluster size variance by its empirical counterpart. However, this is difficult for two reasons: (a) the rareness of the clusters, and (b) the uncertainty on how to group the observed excesses into clusters. For nonparametric estimators of the cluster-size distribution, we refer to (5) and (16).

On the other hand, we can remain in the spirit of the paper and propose a sliding blocks estimator. Recall $N_r(u) = \sum_{i=1}^r I(X_i > u)$ and its compound Poisson limit $N$ in (3.9). Put $\sigma_r^2(u) := \mathrm{var}\big(N_r(u)\big)$. Under an appropriate uniform integrability condition, we have by (3.10), as $r \to \infty$,

$$\sigma_r^2(u_r) = \mathrm{var}\big(N_r(u_r)\big) \to \mathrm{var}(N) = \alpha m_2.$$

Recall $N_{a,b}(u) := \sum_{a<i\leqslant b} I(X_i > u)$. For a threshold $u$ and a block size $r$, we define

$$\bar{N}_{n,r}(u) := \frac{1}{n-r+1}\sum_{i=0}^{n-r} N_{i,i+r}(u),$$

$$\hat{\sigma}_{n,r}^2(u) := \frac{1}{n-2r+1}\sum_{i=0}^{n-r}\big(N_{i,i+r}(u) - \bar{N}_{n,r}(u)\big)^2.$$

We set the denominator equal to $n-2r+1$ in order to reduce the bias of $\hat{\sigma}_{n,r}^2(u)$ from $O(1/k)$ to $O(1/k^2)$ in case $(X_n)_{n\in\mathbb{N}}$ is iid.

The sliding blocks estimator for $c^2 = m_2/m_1^2 - 1 = \theta^2 m_2 - 1 = \theta \operatorname{var}(N)/\tau - 1$ then becomes

$$\hat{c}_{n,r}^2(u) := \frac{\hat{\theta}_{n,r}^{\mathrm{sl}}(u)}{\hat{\tau}_{n,r}(u)} \hat{\sigma}_{n,r}^2(u) - 1.$$

We derive the consistency of $\hat{c}_{n,r}^2(u)$ under a condition on the fourth moment of $N_r(u)$:

$$\mathrm{E}[N_r^4(u_r)] = O(1), \qquad r \to \infty. \tag{5.2}$$

At the price of a longer proof involving a characteristic function argument, condition (4.1) on the moment of order $2p$ (with $1 < p < 2$) would be sufficient as well.

**Proposition 5.1.** *If (3.1), (3.7)–(3.8), (4.1) and (5.2) hold, then* $\hat{c}_{n,r_n}^2(u_{r_n}) \xrightarrow{p} c^2$.

The proof of Proposition 5.1 is given in Appendix C.

## 6. Reducing the bias

Recall $\theta_r(u)$ as in (2.1) and let $\hat{\theta}_{n,r}(u)$ denote either the disjoint or the sliding blocks estimator. The bias of $\hat{\theta}_{n,r}(u)$ can be decomposed into two parts:

$$\mathrm{E}[\hat{\theta}_{n,r}(u)] - \theta = \big(\mathrm{E}[\hat{\theta}_{n,r}(u)] - \theta_r(u)\big) + \big(\theta_r(u) - \theta\big).$$

The component $\theta_r(u) - \theta$ is inherent to the process $(X_n)_{n\in\mathbb{N}}$ itself. For the three examples below it holds that if $(u_r)_{r\in\mathbb{N}}$ is such that $r\bar{F}(u_r) = O(1)$ then $\theta_r(u_r) - \theta = O(1/r)$ as $r \to \infty$.

The presence of the component $\mathrm{E}[\hat{\theta}_{n,r}(u)] - \theta_r(u)$ stems from the fact that the mean of a function of a random vector is in general unequal to this function applied to the mean of the random vector. By a second-order Taylor expansion, it follows that

$$\begin{aligned}\mathrm{E}[\hat{\theta}_{n,r}(u)] - \theta_r(u) \approx{}& \frac{\operatorname{var}\big(\hat{F}_{n,r}(u)\big)}{2\,\mathrm{E}^2[\hat{F}_{n,r}(u)]\,\mathrm{E}[\hat{\tau}_{n,r}(u)]} + \frac{\operatorname{cov}\big(\hat{F}_{n,r}(u), \hat{\tau}_{n,r}(u)\big)}{\mathrm{E}[\hat{F}_{n,r}(u)]\,\mathrm{E}^2[\hat{\tau}_{n,r}(u)]} \\ &- \frac{\log \mathrm{E}[\hat{F}_{n,r}(u)]}{\mathrm{E}^3[\hat{\tau}_{n,r}(u)]} \operatorname{var}\big(\hat{\tau}_{n,r}(u)\big).\end{aligned}$$

By the above expansion and Propositions 3.1 and 3.2, we obtain, as $n \to \infty$,

$$k_n\big(\mathrm{E}[\hat{\theta}_{n,r_n}(u_{r_n})] - \theta_{r_n}(u_{r_n})\big) \to \mu = \begin{cases} \mu_{\mathrm{dj}} & := \theta(2\alpha)^{-1}(e^\alpha - 1) + \alpha^{-1}\theta c^2, \\ \mu_{\mathrm{sl}} & := \theta\alpha^{-2}(e^\alpha - 1 - \alpha) + \alpha^{-1}\theta c^2, \end{cases} \tag{6.1}$$

for the disjoint and the sliding blocks estimator, respectively. Note that $0 \leqslant \mu_{\mathrm{sl}} \leqslant \mu_{\mathrm{dj}}$. If in addition

$$\theta_{r_n}(u_{r_n}) - \theta = o(1/k_n) \tag{6.2}$$

then it follows that, as $n \to \infty$,

$$k_n \big( \mathrm{E}[\hat{\theta}_{n,r}(u)] - \theta \big) \to \mu.$$

Just like the asymptotic variances in Corollary 4.3, the asymptotic biases of the disjoint and sliding blocks estimators in (6.1) are functions of $\theta$, $\alpha = \theta\tau$, and $c^2$. Given consistent estimators of these three quantities, we can estimate $\mu$ and then correct the extremal index estimators by subtracting $\hat{\mu}/k$. Observe that this procedure has to do with the $O(1/k)$ asymptotics of the estimators only, whereas minimization of the asymptotic variance affects the $O(1/\sqrt{k})$ asymptotics.

Note that condition (6.2) is slightly stronger than (4.2). In case $\theta_{r_n}(u_{r_n}) - \theta = o(1/r_n)$, as in the three examples below, (6.2) is equivalent to $k_n = o(r_n)$, that is, $n^{1/2} = o(r_n)$. In contrast, for condition (4.2), the requirement is only that $k_n = o(r_n^2)$, that is, $n^{1/3} = o(r_n)$.

**Example 6.1** (IID sequence)**.** Let $(X_n)_{n \in \mathbb{N}}$ be a sequence of independent random variables with a common, continuous distribution function $F$. Then $\theta = 1$ and

$$r\big(\theta_r(u_r) - 1\big) \to \frac{\tau}{2}, \qquad r \to \infty.$$

**Example 6.2** (Max Auto-Regressive Process)**.** Let $(W_n)_{n \in \mathbb{N}}$ be a sequence of independent, unit-Fréchet distributed random variables. For $0 < \theta \leqslant 1$, let $X_1 = W_1/\theta$ and $X_n = \max\{(1-\theta)X_{n-1}; W_n\}$, $n \geqslant 2$. The extremal index of the sequence is equal to $\theta$ and

$$r\big(\theta_r(u_r) - \theta\big) \to \frac{\tau\theta}{2} + (1-\theta), \qquad r \to \infty.$$

**Example 6.3** (Moving Maximum Process)**.** Let $(W_n)_{n \in \mathbb{N}}$ be a sequence of independent, unit-Fréchet distributed random variables. Let $X_1 = 2W_1$ and $X_n = \max(W_{n-1}, W_n)$. The extremal index of the sequence is equal to $\theta = 1/2$ and

$$r\big(\theta_r(u_r) - \theta\big) \to \frac{\tau}{4}, \qquad r \to \infty.$$

## 7. Numerical examples

### *7.1. Simulation study*

The finite sample properties of the disjoint and sliding blocks estimators for the extremal index are compared in a simulation study. Sequences of length $n = 10\,000$ are simulated from Max Auto-Regressive processes with $\theta = 0.25, 0.5, 0.75$ and 1. For each sequence the estimators $\hat{\theta}^{\mathrm{dj}}_{n,r}(u)$ and $\hat{\theta}^{\mathrm{sl}}_{n,r}(u)$ are computed for five block sizes and two thresholds. The block size is $r = 25, 50, 100, 200$ or 400. The threshold $u$ is the $\lfloor k\tau \rfloor$-th largest order statistic and is defined by either a default value of $\tau = 1$ or the estimate of the optimal value of $\tau$ described in Section 5. The initial estimates of $c^2$ and $\theta$ required in the latter case are based on the threshold when $\tau = 1$. Monte Carlo approximations to the properties of the estimators are computed from 10 000 simulated sequences.

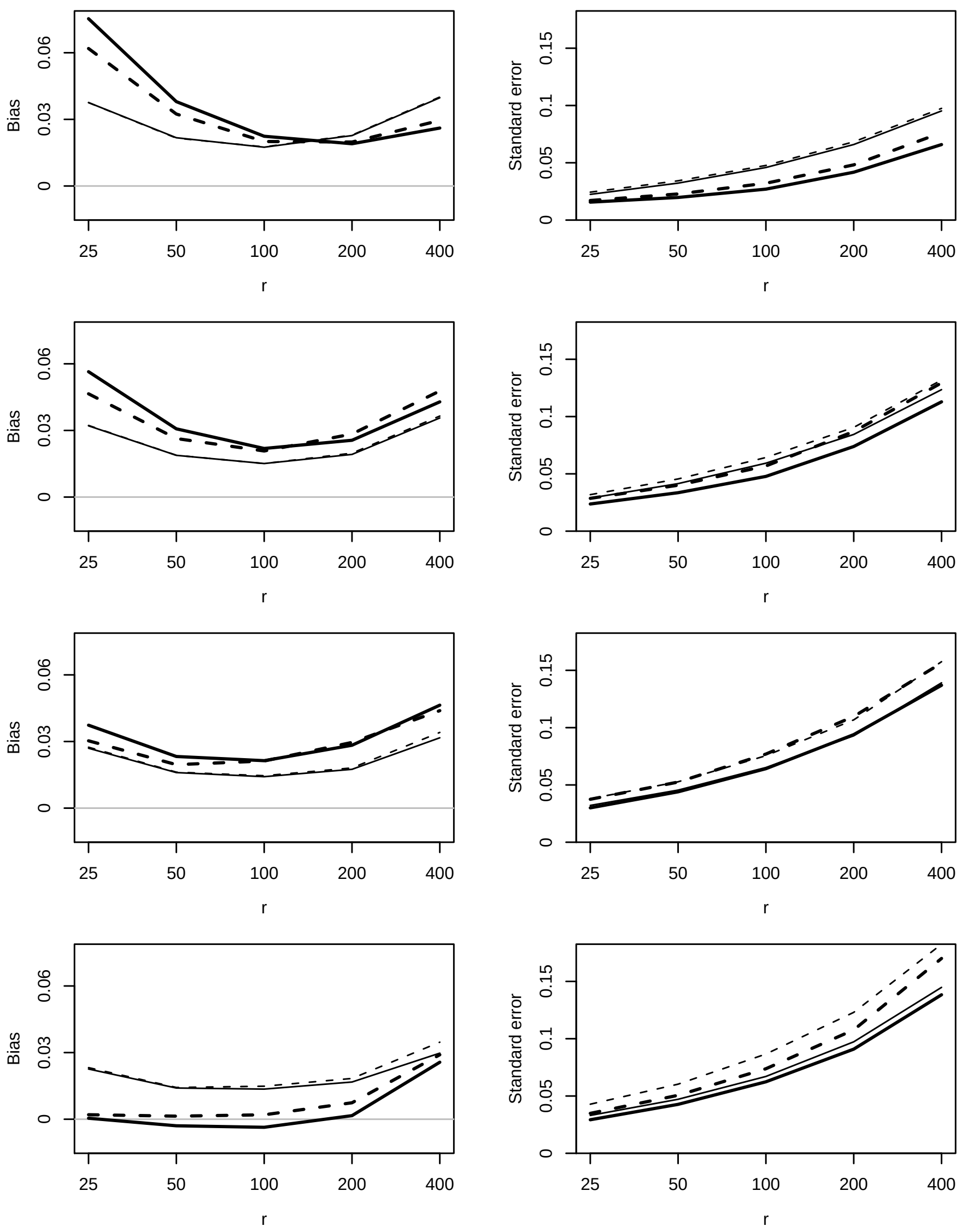

Fig 3. Monte Carlo approximations to the biases (left) and standard errors (right) of the disjoint (dashed) and sliding (solid) blocks estimators for the extremal index with default (thin) and optimal (thick) choices for $\tau$ plotted against block size $r$ on a logarithmic scale. Data are simulated from Max Auto-Regressive processes with $\theta = 0.25$, 0.5, 0.75 and 1 (top to bottom).

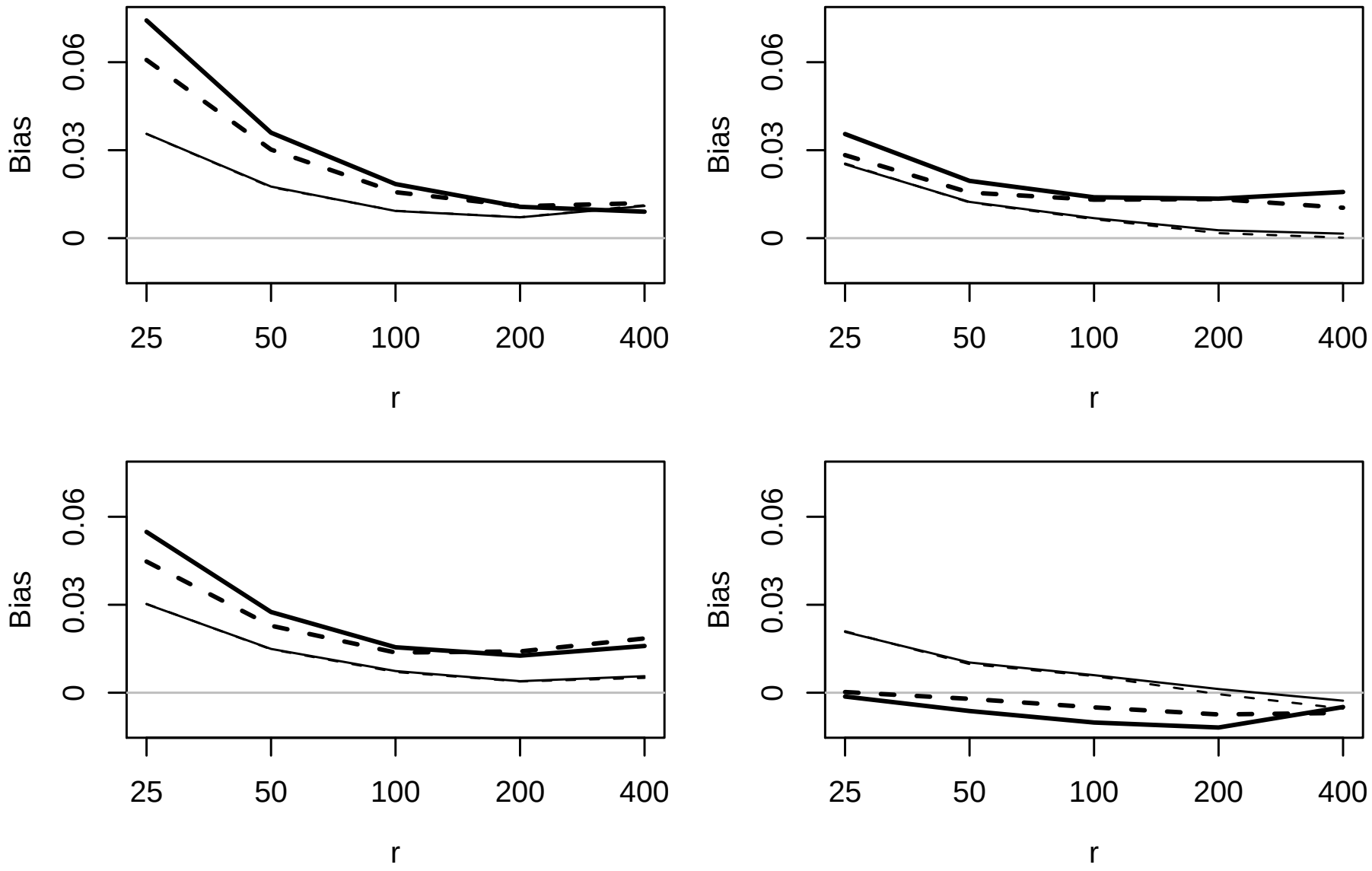

Fig 4. Monte Carlo approximations to the biases of the bias-corrected disjoint (dashed) and sliding (solid) blocks estimators for the extremal index with default (thin) and optimal (thick) choices for $\tau$ plotted against block size $r$ on a logarithmic scale. Data are simulated from Max Auto-Regressive processes with $\theta = 0.25$ (top-left), 0.5 (bottom-left), 0.75 (top-right) and 1 (bottom-right).

Figure 3 shows the biases and standard errors of the estimators. Biases tend to be positive and smallest at intermediate block sizes while variances increase with block size. Sliding blocks always yield lower standard errors than disjoint blocks. There is also evidence that sliding blocks yield larger biases than disjoint blocks when $r$ is small and smaller biases when $r$ is large. Optimizing $\tau$ tends to yield lower variances than the default $\tau = 1$, but also larger biases when $\theta < 1$. This is explained by the fact that the estimated values for the optimal $\tau$ tend to exceed 1 except when $\theta = 1$. Example 6.2 suggests that increasing $\tau$ increases the bias.

The effect of the bias correction described in Section 6 is shown in Figure 4. There is little improvement for small block sizes, but biases are reduced significantly and stabilized for larger block sizes. The impact on the standard errors is negligible (not shown).

The positive biases of $\hat{\theta}^{\mathrm{dj}}_{n,r}(u)$ and $\hat{\theta}^{\mathrm{sl}}_{n,r}(u)$ can lead to poor coverage properties (not shown) of confidence intervals for $\theta$ based on the asymptotic Normal distribution of Section 4. Lower and upper confidence limits tend to be too high when $r$ is small but coverage improves when $r$ is large. Coverage is also affected

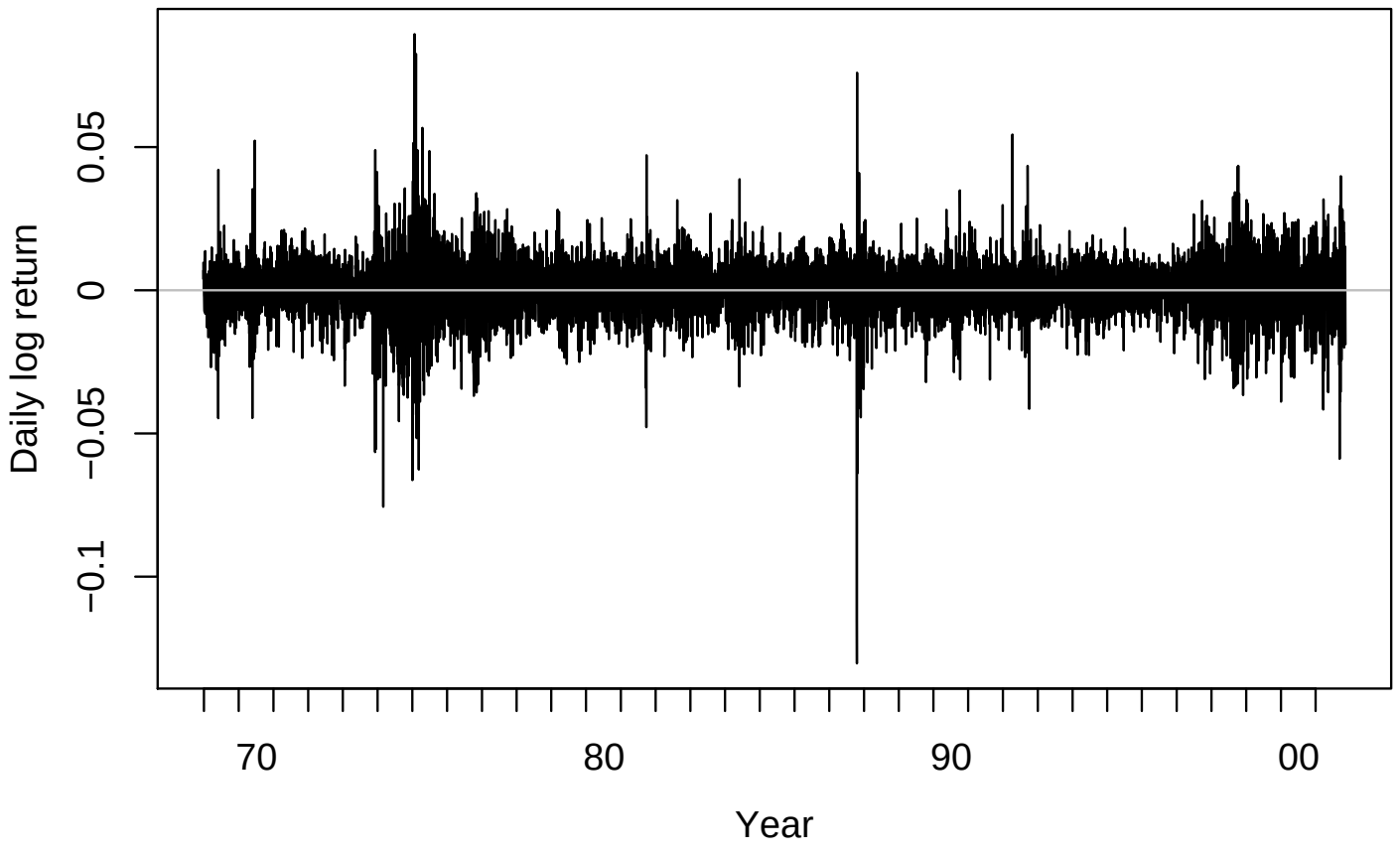

FIG 5. Daily log returns of the FTSE100 index from 25 December 1968 to 12 November 2001.

by underestimation of standard errors when $\theta < 1$ (not shown).

These simulations were repeated for the doubly stochastic process of (17) and for ARCH(1) processes; see (3) and **(author?)** (4, Chapter 8). Results for the doubly stochastic process were very similar to those reported above. Results for ARCH(1) processes were also similar but the improvement in variance afforded by sliding blocks was less clear. Qualitatively similar results were found when the simulations were repeated with $n = 1000$ and $r = 5, 10, 20, 40$ or 80.

### *7.2. Case study*

The extremal index is now estimated for a financial time series: daily log returns of the FTSE100 index between 25 December 1968 and 12 November 2001. This series was analysed previously by (10) and is plotted in Figure 5; the data were kindly passed on by Jonathan Tawn. Clusters of large, negative returns can be financially damaging so estimates of the extremal index for the negated series are plotted against block size in Figure 6. Two sliding blocks estimators are compared: both employ the bias correction but one uses the default value $\tau = 1$ while the other uses estimated optimal values $\tau = \hat{\tau}_{\text{opt}}$. Thresholds are the $\lfloor k\tau \rfloor$-th largest order statistics so that the proportion of data exceeding the threshold for block size $r$ is $\tau/r$. The lower horizontal axis in Figure 6 is therefore a transformation of the threshold used when $\tau = 1$, and coincides with the scale used by (10). The upper horizontal axis represents the same transformation of the threshold when $\tau = \hat{\tau}_{\text{opt}}$. These latter thresholds are lower because $\hat{\tau}_{\text{opt}} \approx 5$ for all but the smallest block sizes.

The point estimates from the two sliding blocks estimators are similar and both stablize near $\theta = 1/3$. Estimates from the intervals estimator of (6) are also shown and differ slightly but are consistent with an extremal index of one-third once sampling variation is taken into account. However, these values are approximately half those obtained by (10) with a two-thresholds estimator. The confidence intervals in Figure 6 are computed using the estimated standard errors for the sliding blocks estimators with no bias correction. The confidence intervals when $\tau = \hat{\tau}_{\mathrm{opt}}$ are often much narrower than when $\tau = 1$ owing to the lower thresholds mentioned above.

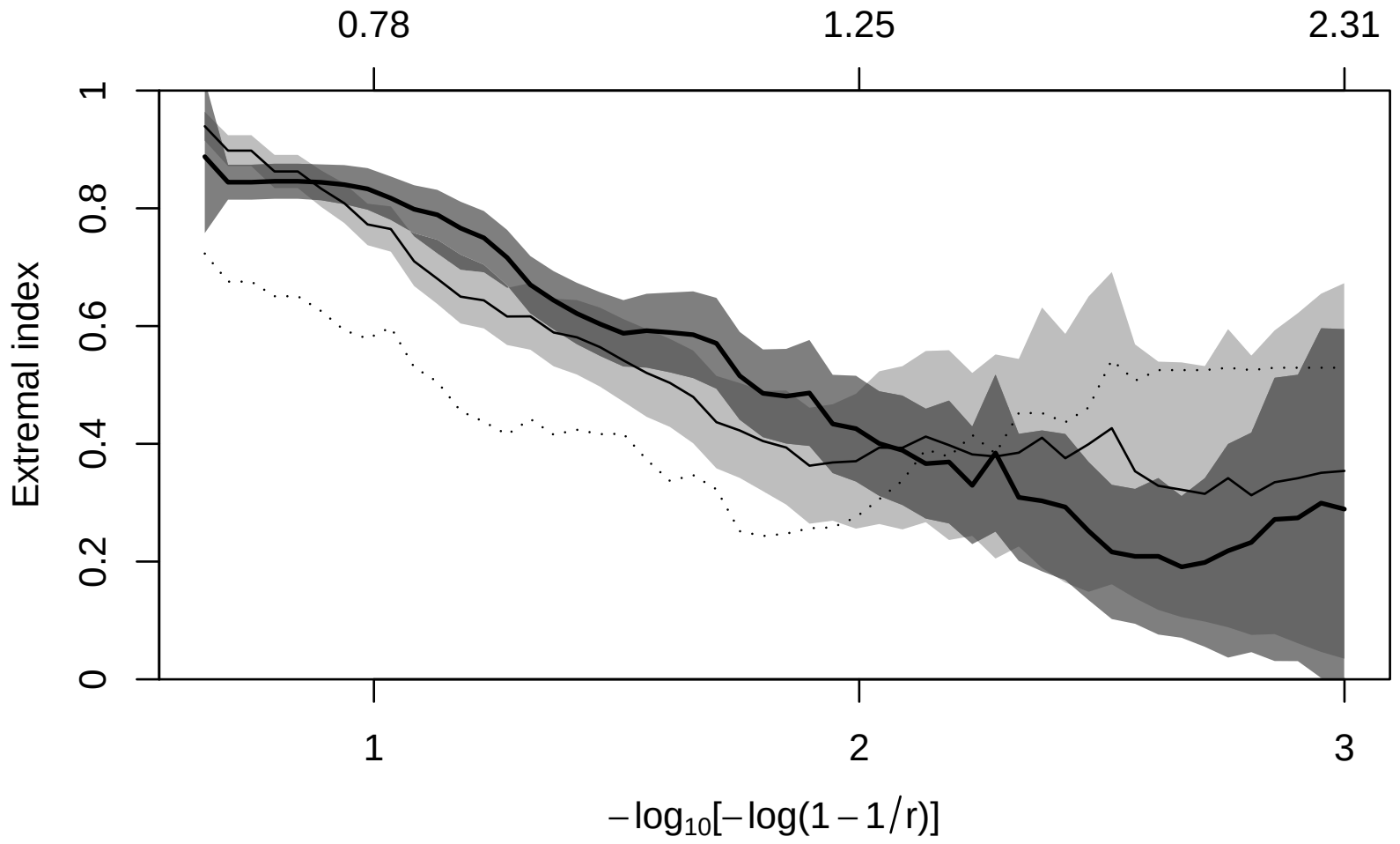

Fig 6. Sliding blocks estimates (solid) of the extremal index for the negative daily log returns plotted against block size $r$ on a log-log scale. The estimates use default ($\tau = 1$; thin) and optimal ($\tau = \hat{\tau}_{\mathrm{opt}}$; thick) choices for $\tau$. Intervals estimates (dotted) are based on the threshold obtained when $\tau = 1$. Shading indicates pointwise 90% confidence intervals based on the asymptotic Normal approximation for the sliding blocks estimates. The top axis represents $-\log_{10}[-\log(1 - \hat{\tau}_{\mathrm{opt}}/r)]$.

## Appendix A: Proofs for Section 3

### A.1. Proof of Proposition 3.1

*Asymptotic variance of* $\hat{F}^{\mathrm{dj}}_{n,r_n}(u_{r_n})$. By stationarity,

$$\begin{aligned}
\operatorname{var}\big(\hat{F}^{\mathrm{dj}}_{n,r}(u)\big) &= \frac{1}{k} F_r(u)\{1 - F_r(u)\} \\
&\quad + \frac{2}{k^2} \sum\sum_{1 \leqslant i < j \leqslant k} \{\Pr(M_{(i-1)r,ir} \leqslant u, M_{(j-1)r,jr} \leqslant u) - F_r^2(u)\} \\
&= \frac{1}{k} F_r(u)\{1 - F_r(u)\} + \frac{2(k-1)}{k^2}\{F_{2r}(u) - F_r^2(u)\} \\
&\quad + \frac{2}{k^2} \sum_{s=2}^{k-1} (k-s)\{\Pr(M_{0,r} \leqslant u, M_{sr,(s+1)r} \leqslant u) - F_r^2(u)\}.
\end{aligned}$$

By definition of the extremal index, $F_{r_n}(u_{r_n}) \to e^{-\tau\theta}$ and $F_{2r_n}(u_{r_n}) \to e^{-2\tau\theta}$ as $n \to \infty$. Hence (3.2) will follow if we can show that

$$\frac{1}{k_n} \sum_{s=2}^{k_n-1} (k_n - s)\{\Pr(M_{0,r_n} \leqslant u_{r_n}, M_{sr_n,(s+1)r_n} \leqslant u_{r_n}) - F_{r_n}^2(u_{r_n})\} \to 0$$

as $n \to \infty$. But since the maximal correlation coefficient $\rho_{n,l}(u)$ is decreasing in $l$, the left-hand of the previous display is bounded in absolute value by

$$\sum_{s=2}^{k_n-1} \rho_{n,(s-1)r_n}(u_{r_n}) \leqslant \frac{1}{r_n} \sum_{l=1}^{n} \rho_{n,l}(u_{r_n}).$$

By hypothesis, this converges to zero as $n \to \infty$. □

*Asymptotic variance of* $\hat{F}^{\mathrm{sl}}_{n,r_n}(u_{r_n})$. By stationarity,

$$\begin{aligned}
\operatorname{var}\big(\hat{F}^{\mathrm{sl}}_{n,r}(u)\big) &= \frac{1}{n-r+1} F_r(u)\{1 - F_r(u)\} \\
&\quad + \frac{2}{(n-r+1)^2} \sum\sum_{0 \leqslant i < j \leqslant n-r} \operatorname{cov}\big(I(M_{i,i+r} \leqslant u), I(M_{j,j+r}) \leqslant u\big).
\end{aligned}$$

The sum on the right-hand side of the previous display can be written as

$$\begin{aligned}
&\frac{2}{(n-r+1)^2} \sum_{s=1}^{n-r} (n-r+1-s) \operatorname{cov}\big(I(M_{0,r} \leqslant u), I(M_{s,s+r} \leqslant u)\big) \\
&= \frac{2r}{n-r+1} \frac{1}{r} \sum_{s=1}^{r} \left(1 - \frac{s}{n-r+1}\right) \{F_{s+r}(u) - F_r(u)^2\} \\
&\quad + \frac{2r}{n-r+1} \frac{1}{r} \sum_{s=r+1}^{n-r} \left(1 - \frac{s}{n-r+1}\right) \operatorname{cov}\big(I(M_{0,r} \leqslant u), I(M_{s,s+r} \leqslant u)\big).
\end{aligned}$$

By dominated convergence, as $n \to \infty$,

$$\frac{1}{r_n}\sum_{s=1}^{r_n}\left(1-\frac{s}{n-r_n+1}\right)\{F_{s+r_n}(u_{r_n})-F_{r_n}(u_{r_n})^2\}$$
$$\to \int_0^1 (e^{-(u+1)\tau\theta}-e^{-2\tau\theta})\,\mathrm{d}u = e^{-\tau\theta}\frac{1-e^{-\tau\theta}}{\tau\theta}-e^{-2\tau\theta}.$$

Hence (3.3) will follow if we can show that, as $n \to \infty$,

$$\frac{1}{r_n}\sum_{s=r_n+1}^{n-r_n}\left(1-\frac{s}{n-r_n+1}\right)\mathrm{cov}\big(I(M_{0,r_n}\leqslant u_{r_n}), I(M_{s,s+r_n}\leqslant u_{r_n})\big)\to 0.$$

But this follows from the assumption that $r_n^{-1}\sum_{l=1}^n \rho_{n,l}(u_{r_n}) \to 0$ as $n \to \infty$. □

*Asymptotic covariance of $\hat{F}^{\mathrm{dj}}_{n,r_n}(u_{r_n})$ and $\hat{F}^{\mathrm{sl}}_{n,r_n}(u_{r_n})$.* We have

$$\begin{aligned} &k\,\mathrm{cov}\big(\hat{F}^{\mathrm{dj}}_{n,r}(u),\hat{F}^{\mathrm{sl}}_{n,r}(u)\big)\\ &=\mathrm{cov}\big(I(M_{0,r}\leqslant u),\hat{F}^{\mathrm{sl}}_{n,r}(u)\big)+\mathrm{cov}\big(I(M_{(k-1)r,kr}\leqslant u),\hat{F}^{\mathrm{sl}}_{n,r}(u)\big)\\ &\quad+\frac{1}{n-r+1}\sum_{i=2}^{k-1}\sum_{j=0}^{n-r}\mathrm{cov}\big(I(M_{(i-1)r,ir}\leqslant u),I(M_{j,j+r}\leqslant u)\big).\end{aligned}$$

The first two terms on the right-hand side are bounded by $\{\mathrm{var}(\hat{F}^{\mathrm{sl}}_{n,r}(u))\}^{1/2}$; in view of (3.3), they will not contribute to the limit. The final term on the right-hand side of the previous display can be decomposed into two pieces, $I+II$ say, according to whether $(i-2)r \leqslant j < ir$ or not. For the first term, the union of the intervals of integers $\{(i-1)r+1,\ldots,ir\}$ and $\{j+1,\ldots,j+r\}$ is again an interval of integers; by stationarity,

$$\begin{aligned} I &= \frac{k-2}{n-r+1}\sum_{j=0}^{2r-1}\mathrm{cov}\big(I(M_{r,2r}\leqslant u),I(M_{j,j+r}\leqslant u)\big)\\ &= \frac{r(k-2)}{n-r+1}\frac{1}{r}\sum_{j=0}^{2r-1}\{F_{(2r-j)\vee j}(u)-F_r^2(u)\}.\end{aligned}$$

Adding the subscript $n$ to indicate the dependence on $n$, we get

$$\lim_{n\to\infty} I_n = \int_0^2 (e^{-\{(2-u)\vee u\}\alpha}-e^{-2\alpha})\,\mathrm{d}u = 2\left(e^{-\alpha}\frac{1-e^{-\alpha}}{\alpha}-e^{-2\alpha}\right).$$

The second term can be bounded as follows

$$|II| \leqslant \frac{2(k-2)}{n-r+1}\sum_{l=1}^{n}\rho_{n,l}(u).$$

Since $r_n^{-1}\sum_{l=1}^n \rho_{n,l}(u_{r_n}) \to 0$ by assumption, we get $II_n \to 0$ as $n \to \infty$. □

### A.2. Proof of Proposition 3.2

*Asymptotic covariance of* $\hat{F}^{\mathrm{dj}}_{n,r_n}(u_{r_n})$ *and* $\hat{\tau}_{n,r_n}(u_{r_n})$. We have

$$k \operatorname{cov}\bigl(\hat{F}^{\mathrm{dj}}_{n,r}(u), \hat{\tau}_{n,r}(u)\bigr) = -\frac{1}{k} \sum_{i=1}^{k} \sum_{j=1}^{kr} \operatorname{cov}\bigl(I(M_{(i-1)r,ir} \leqslant u), I(X_j \leqslant u)\bigr).$$

We split the sum according into two pieces, $I + II$ say, according to whether $(i-1)r < j \leqslant r$ or not. By stationarity, the first term is equal to

$$I = -\sum_{j=1}^{r} \operatorname{cov}\bigl(I(M_{0,r} \leqslant u), I(X_j \leqslant j)\bigr) = -r\{1 - F(u)\}F_r(u). \tag{A.1}$$

Adding subscripts $n$ to indicate the dependence on $n$, we get $I_n \to -\tau e^{-\alpha}$ as $n \to \infty$. The second term, $II$, can be bounded in absolute value as follows:

$$\begin{aligned} |II| &\leqslant 2\bigl\{\operatorname{var}\bigl(I(M_{0,r} \leqslant u)\bigr) \operatorname{var}\bigl(I(X_1 \leqslant u)\bigr)\bigr\}^{1/2} \sum_{l=1}^{n} \rho_{n,l}(u) \\ &\leqslant 2\{1 - F(u)\}^{1/2} \sum_{l=1}^{n} \rho_{n,l}(u). \end{aligned} \tag{A.2}$$

Assumptions (3.1) and (3.5) now imply that $II_n \to 0$ as $n \to \infty$. □

*Asymptotic covariance of* $\hat{F}^{\mathrm{sl}}_{n,r_n}(u_{r_n})$ *and* $\hat{\tau}_{n,r_n}(u_{r_n})$. This time, we have

$$k \operatorname{cov}\bigl(\hat{F}^{\mathrm{sl}}_{n,r}(u), \hat{\tau}_{n,r}(u)\bigr) = -\frac{1}{n-r+1} \sum_{i=0}^{n-r} \sum_{j=1}^{kr} \operatorname{cov}\bigl(I(M_{i,i+r} \leqslant u), I(X_j \leqslant u)\bigr).$$

We split the sum according into two pieces, $I + II$ say, according to whether $i < j \leqslant i+r$ or not. The first term, $I$ is the same as in (A.1) and so gives rise to the same limit. The second term, $II$, admits the same bounds as in (A.2) and hence is asymptotically negligible. □

### A.3. Proof of Proposition 3.3

Recall $N_s(u) = \sum_{i=1}^{s} I(X_i > u)$. We have

$$k \operatorname{var}\bigl(\hat{\tau}_{n,r}(u)\bigr) = \frac{1}{k} \operatorname{var}\bigl(N_{rk}(u)\bigr).$$

For integer $0 \leqslant a \leqslant b$, put $N_{a,b}(u) := \sum_{a<i\leqslant b} I(X_i > u)$, the sum being zero if $a = b$. Fix integer $1 \leqslant l < s < n$ and write $m := \lfloor rk/s \rfloor$. We have

$$\begin{aligned} N_{r_n k_n}(u_{r_n}) = \sum_{i=1}^{m} N_{(i-1)s_n, is_n - l_n}(u_{r_n}) + \sum_{i=1}^{m} N_{is_n - l_n, is_n}(u_{r_n}) \\ + N_{m_n s_n, r_n k_n}(u_{r_n}) =: A_n + B_n + C_n. \end{aligned}$$

By the Cauchy-Schwarz inequality, it is sufficient to show that, as $n \to \infty$,

$$(1/k_n)\operatorname{var}(A_n) \to \alpha m_2, \qquad (1/k_n)\operatorname{var}(B_n) \to 0, \qquad (1/k_n)\operatorname{var}(C_n) \to 0.$$

Before we treat these three terms, it is useful to note that the assumptions imply that

$$\operatorname{var}\bigl(N_{l_n}(u_{r_n})\bigr) = o\bigl(s_n \Pr(X_1 > u_{r_n})\bigr), \qquad n \to \infty, \tag{A.3}$$

which in turns implies

$$(1/k_n) m_n \operatorname{var}\bigl(N_{l_n}(u_{r_n})\bigr) \to 0, \qquad n \to \infty. \tag{A.4}$$

[Proof of (A.3): $\operatorname{var}\bigl(N_{l_n}(u_{r_n})\bigr) \leqslant \mathrm{E}[N_{l_n}^2(u_{r_n})] \leqslant s_n \Pr(X_1 > u_{r_n}) \mathrm{E}[N_{l_n}^2(u_{r_n}) \mid M_{s_n} > u_{r_n}]$ and the final expectation tends to zero by uniform integrability and $\Pr\{N_{l_n}(u_{r_n}) > 0 \mid M_{s_n} > u_{r_n}\} \leqslant l_n \Pr(X_1 > u_{r_n}) / \Pr(M_{s_n} > u_{r_n}) \sim (l_n/s_n)\theta^{-1} \to 0$.]

*The term $A_n$.* By stationarity, $(1/k_n)\operatorname{var}(A_n) = I_n + II_n$ with

$$\begin{aligned} I_n &:= (1/k_n) m_n \operatorname{var}\bigl(N_{s_n - l_n}(u_{r_n})\bigr), \\ II_n &:= 2(1/k_n) \sum_{1 \leqslant i < j \leqslant m} \sum \operatorname{cov}\bigl(N_{(i-1)s_n, is_n - l_n}(u_{r_n}), N_{(j-1)s_n, js_n - l_n}(u_{r_n})\bigr). \end{aligned}$$

We first treat $I_n$. We have

$$\begin{aligned} &\operatorname{var}\bigl(N_{s_n}(u_{r_n})\bigr) \\ &\qquad = \mathrm{E}[N_{s_n}^2(u_{r_n}) \mid M_{s_n} > u_{r_n}] \Pr(M_{s_n} > u_{r_n}) - \{s_n \Pr(X_1 > u_{r_n})\}^2. \end{aligned}$$

Since $m_n \sim n/s_n$ and $\Pr(M_{s_n} > u_{r_n}) \sim s_n \Pr(X_1 > u_{r_n})\theta$ as $n \to \infty$, we find

$$(1/k_n) m_n \operatorname{var}\bigl(N_{s_n}(u_{r_n})\bigr) \to \theta\tau m_2, \qquad n \to \infty.$$

Since $N_{s_n}(u_{r_n}) = N_{s_n - l_n}(u_{r_n}) + N_{s_n - l_n, s_n}(u_{r_n})$ and since $N_{s_n - l_n, s_n}(u_{r_n})$ and $N_{l_n}(u_{r_n})$ have the same distribution, the previous display and (A.4) imply that $I_n \to \theta\tau m_2$ as $n \to \infty$.

Next we treat $II_n$. We have

$$|II_n| \leqslant 2(1/k_n) m_n^2 \operatorname{var}\bigl(N_{s_n - l_n}(u_{r_n})\bigr) \rho_{n, l_n}(u_{r_n}).$$

In view of what we obtained for $I_n$ and since $m_n \rho_{n, l_n}(u_{r_n}) \to 0$, we conclude that $II_n \to 0$ as $n \to \infty$.

*The term $B_n$.* By stationarity

$$\begin{aligned} &(1/k_n)\operatorname{var}(B_n) \\ &\qquad \leqslant (1/k_n) m_n \operatorname{var}\bigl(N_{l_n}(u_{r_n})\bigr) + 2(1/k_n) m_n^2 \operatorname{var}\bigl(N_{l_n}(u_{r_n})\bigr) \rho_{n, l_n}(u_{r_n}). \end{aligned}$$

By (A.4) and since $m_n \rho_{n, l_n}(u_{r_n}) \to 0$, we obtain that $(1/k_n)\operatorname{var}(B_n) \to 0$ as $n \to \infty$.

*The term $C_n$.* By stationarity,

$$\begin{aligned}(1/k_n)\operatorname{var}(C_n) &= (1/k_n)\operatorname{var}\bigl(N_{n-m_ns_n}(u_{r_n})\bigr)\\ &\leqslant (1/k_n)\operatorname{E}[N^2_{n-m_ns_n}(u_{r_n})] \leqslant (1/k_n)\operatorname{E}[N^2_{s_n}(u_{r_n})]\\ &= (1/k_n)\Pr(M_{s_n} > u_{r_n})\operatorname{E}[N^2_{s_n}(u_{r_n}) \mid M_{s_n} > u_{r_n}].\end{aligned}$$

By assumption, the limit as $n \to \infty$ is zero. □

## Appendix B: Proofs for Section 4

### B.1. *Proof of Theorem 4.1*

By Propositions 3.1 and 3.3 in combination with Tchebychev's inequality, it is not difficult to see that as $n \to \infty$,

$$\bigl|\hat{F}^{\mathrm{dj}}_{n,r_n}(u_{r_n}) - F_{r_n}(u_{r_n})\bigr| + \bigl|\hat{F}^{\mathrm{sl}}_{n,r_n}(u_{r_n}) - F_{r_n}(u_{r_n})\bigr| + \bigl|\hat{\tau}_{n,r_n}(u_{r_n}) - \tau_{r_n}(u_{r_n})\bigr| \xrightarrow{p} 0.$$

By definition of the extremal index, $F_{r_n}(u_{r_n}) \to e^{-\theta\tau}$ and $\tau_{r_n}(u_{r_n}) \to \tau$ as $n \to \infty$. The result follows by continuity of $(x, y) \mapsto -\log(x)/y$ on $(0,\infty) \times (0,\infty)$. □

### B.2. *Proof of Theorem 4.2*

Write $\boldsymbol{Z}_n = (Z_{n,1}, Z_{n,2}, Z_{n,3})^\top$ with

$$\begin{aligned}&Z_{n,1} = \sqrt{k_n}\bigl(\hat{F}^{\mathrm{dj}}_{n,r_n}(u_{r_n}) - F_{r_n}(u_{r_n})\bigr), \quad Z_{n,3} = \sqrt{k}\bigl(\hat{\tau}_{n,r_n}(u_{r_n}) - \tau_{r_n}(u_{r_n})\bigr),\\ &Z_{n,2} = \sqrt{k_n}\bigl(\hat{F}^{\mathrm{sl}}_{n,r_n}(u_{r_n}) - F_{r_n}(u_{r_n})\bigr).\end{aligned}$$

By the Cramér-Wold device, it is sufficient to show that for $\boldsymbol{a} = (a_1, a_2, a_3)^\top \in \mathbb{R}^3$,

$$a_1 Z_{n,1} + a_2 Z_{n,2} + a_3 Z_{n,3} \xrightarrow{d} N(0, \boldsymbol{a}^\top \Sigma \boldsymbol{a}), \qquad n \to \infty.$$

Note that

$$Z_{n,1} = \frac{1}{\sqrt{k_n}} \sum_{i=1}^{k_n} \bar{I}^{\mathrm{dj}}_{i,r_n}, \qquad Z_{n,3} = \frac{1}{\sqrt{k_n}} \sum_{i=1}^{k_n} \bar{N}_{i,r_n},$$

$$Z_{n,2} = \frac{\sqrt{k_n}}{n - r_n + 1} \sum_{t=0}^{n-r_n} \bar{I}^{\mathrm{sl}}_{t,r_n},$$

where for $i \in \{1, \ldots, k\}$ and $t \in \{0, \ldots, n-r\}$,

$$\begin{aligned}&\bar{I}^{\mathrm{dj}}_{i,r} = I(M_{(i-1)r,ir} \leqslant u_r) - F_r(u_r), \quad \bar{N}_{i,r} = \sum_{t=(i-1)r+1}^{ir} \bigl(I(X_t > u_r) - \bar{F}(u_r)\bigr),\\ &\bar{I}^{\mathrm{sl}}_{t,r} = I(M_{t,t+r} \leqslant u_r) - F_r(u_r).\end{aligned}$$

The idea of the proof is as follows: By clipping out certain terms in the definitions of $Z_{n,j}$, the latter can be viewed upon as sums of approximately independent random variables. Asymptotic normality then follows from the Lindeberg-Feller central limit theorem for triangular arrays.

Let $k^* \in \{1, \ldots, k-1\}$. Construct a partition of $\{1, \ldots, k\}$ into subsets of size $k^*$, with two adjacent such subsets separated by a singleton. The number of subsets of size $k^*$ that can be formed in this way is $q = \lfloor (k+1)/(k^*+1) \rfloor$. We have

$$\frac{1}{\sqrt{k}} \sum_{i=1}^{k} \bar{I}_{i,r}^{\mathrm{dj}} = \frac{1}{\sqrt{k}} \sum_{j=1}^{q} \sum_{i=(j-1)(k^*+1)+1}^{j(k^*+1)-1} \bar{I}_{i,r_n}^{\mathrm{dj}} + \frac{1}{\sqrt{k}} \sum_{j=1}^{q-1} \bar{I}_{j(k^*+1),r_n}^{\mathrm{dj}} + \frac{1}{\sqrt{k}} \sum_{i=q(k^*+1)}^{k} \bar{I}_{i,r_n}^{\mathrm{dj}}.$$

Let $k^* = k_n^*$ be such that $k_n^* \to \infty$ but $k_n^* = o(k_n)$ as $n \to \infty$. The final two terms on the right-hand side of the previous display are negligible as their variances tend to zero: the variance of the second term on the right-hand side is of the order

$$O\left( \frac{q_n}{k_n} \operatorname{var}(\bar{I}_{1,r_n}^{\mathrm{dj}}) + 2 \frac{q_n^2}{k_n} \rho_{n,r_n}(u_{r_n}) \operatorname{var}(\bar{I}_{1,r_n}^{\mathrm{dj}}) \right) = o(1),$$

where we used $q_n = o(k_n)$, $\operatorname{var}(\bar{I}_{1,r_n}^{\mathrm{dj}}) = O(1)$, and

$$\rho_{n,r_n}(u_{r_n}) = o(1/k_n) \tag{B.1}$$

[Proof of (B.1): by (3.7) $\rho_{n,r_n}(u_{r_n}) \leqslant \rho_{n,l_n}(u_{r_n}) = o(s_n/n)$ and $s_n = o(r_n)$]; similarly, the variance of the third term is of the order

$$O\left( 2 \frac{k_n^*}{k_n} \operatorname{var}(\bar{I}_{1,r_n}^{\mathrm{dj}}) + 2 \frac{(k_n^*)^2}{k_n} \rho_{n,r_n}(u_{r_n}) \operatorname{var}(\bar{I}_{1,r_n}^{\mathrm{dj}}) \right) = o(1).$$

As a consequence,

$$Z_{n,1} = \frac{1}{\sqrt{k_n}} \sum_{j=1}^{q_n} \sum_{i=(j-1)(k_n^*+1)+1}^{j(k_n^*+1)-1} \bar{I}_{i,r_n}^{\mathrm{dj}} + o_p(1). \tag{B.2}$$

In a completely similar way, just replacing $\bar{I}_{i,r}^{\mathrm{dj}}$ by $\bar{N}_{i,r}$, we can also show that

$$Z_{n,3} = \frac{1}{\sqrt{k_n}} \sum_{j=1}^{q_n} \sum_{i=(j-1)(k_n^*+1)+1}^{j(k_n^*+1)-1} \bar{N}_{i,r_n} + o_p(1); \tag{B.3}$$

a crucial element here is that $\operatorname{var}(\bar{N}_{1,r_n}) = O(1)$, which follows from (4.1).

Next, construct a partition $\{1, \ldots, n\}$ into $k$ blocks of size $r$ and in case $kr < n$ a final block of length $n - kr$. Form subsamples by taking unions over $k^*$ consecutive blocks of size $r$, two consecutive subsamples being separated by a single block of size $r$. The number of subsamples that can formed in this way is again $q = \lfloor (k+1)/(k^*+1) \rfloor$, the $j$th subsample being

$$\begin{aligned}\mathcal{S}_{j,k^*,r} &= \bigcup_{i=(j-1)(k^*+1)+1}^{j(k^*+1)-1} \{(i-1)r+1, \ldots, ir\} \\ &= \{(j-1)(k^*+1)r+1, \ldots, j(k^*+1)r - r\}.\end{aligned}$$

In the definition of the sliding-blocks estimator, retain only those $t \in \{0, \ldots, n-r\}$ such that the (sliding) block $\{t+1, \ldots, t+r\}$ is contained entirely in one of the subsamples $\mathcal{S}_{j,k^*,r}$. In other words, discard those $t$ such that $\{t+1, \ldots, t+r\}$ has a nonempty intersection with one of the $q-1$ blocks of size $r$ separating two consecutive subsamples or with the remaining part of the sample after the final subsample $\mathcal{S}_{q,k^*,r}$. The values of $t$ to be retained are then given as follows: for $j \in \{1, \ldots, q\}$,

$$\{t+1, \ldots, t+r\} \subset \mathcal{S}_{j,k^*,r} \text{ if and only if } (j-1)(k^*+1)r \leqslant t \leqslant j(k^*+1)r - 2.$$

We find

$$\begin{aligned}\frac{\sqrt{k}}{n} \sum_{t=0}^{n-r} \bar{I}^{\mathrm{sl}}_{t,r} &= \frac{\sqrt{k}}{n} \sum_{j=1}^{q} \sum_{t=(j-1)(k^*+1)r}^{j(k^*+1)r-2r} \bar{I}^{\mathrm{sl}}_{t,r} \\ &\quad + \frac{\sqrt{k}}{n} \sum_{j=1}^{q-1} \sum_{t=j(k^*+1)r-2r+1}^{j(k^*+1)r-1} \bar{I}^{\mathrm{sl}}_{t,r} + \frac{\sqrt{k}}{n} \sum_{t=q(k^*+1)r-2r+1}^{n-r} \bar{I}^{\mathrm{sl}}_{t,r}.\end{aligned}$$

Again, let $k^* = k_n^*$ be such that $k_n^* \to \infty$ and $k_n^* = o(k_n)$ as $n \to \infty$. Asymptotically, the variances of the final two sums tend to zero: the variance of the second sum on the right-hand side is of the order

$$O\left(\frac{k_n}{n^2}\{q_n + q_n^2 \rho_{n,r_n}(u_{r_n})\} r_n^2\right) = O\left(\frac{q_n}{k_n} + \frac{q_n^2}{k_n} \rho_{n,r_n}(u_{r_n})\right) = o(1);$$

and since the number of terms in the third sum on the right-hand side is not larger than $2k^*r$, the variance of that sum is of the order the variance of that sum is of the order

$$O\left(\frac{k_n}{n^2} k_n^* r_n^2\right) = O\left(\frac{k_n^*}{k_n}\right) = o(1).$$

As a consequence,

$$Z_{n,2} = \frac{\sqrt{k_n}}{n - r_n + 1} \sum_{j=1}^{q_n} \sum_{t=(j-1)(k_n^*+1)r_n}^{j(k_n^*+1)r_n - 2r_n} \bar{I}^{\mathrm{sl}}_{t,r_n} + o_p(1). \tag{B.4}$$

Equations (B.2), (B.3) and (B.4) can be summarised as

$$Z_{n,v} = \frac{1}{\sqrt{q_n}} \sum_{j=1}^{q_n} \xi_{n,j,v} + o_p(1), \qquad v \in \{1,2,3\},$$

with for $j \in \{1,\ldots,q_n\}$,

$$\xi_{n,j,1} = \sqrt{q_n/k_n} \sum_{i=(j-1)(k_n^*+1)+1}^{j(k_n^*+1)-1} \bar{I}^{\mathrm{dj}}_{i,r_n},$$

$$\xi_{n,j,2} = \frac{\sqrt{q_n k_n}}{n-r_n+1} \sum_{t=(j-1)(k_n^*+1)r_n}^{j(k_n^*+1)r_n-2r_n} \bar{I}^{\mathrm{sl}}_{t,r_n},$$

$$\xi_{n,j,3} = \sqrt{q_n/k_n} \sum_{i=(j-1)(k_n^*+1)+1}^{j(k_n^*+1)-1} \bar{N}_{i,r_n}.$$

As a consequence,

$$a_1 Z_{n,1} + a_2 Z_{n,2} + a_3 Z_{n,3} = \frac{1}{\sqrt{q_n}} \sum_{j=1}^{q_n} \xi_{n,j} + o_p(1)$$

with $\xi_{n,j} = a_1\xi_{n,j,1} + a_2\xi_{n,j,2} + a_3\xi_{n,j,3}$. Note that $\xi_{n,j,v}$ is measurable with respect to the $\sigma$-field generated by the events $\{X_t \leqslant u_r\}$ with $t$ ranging over the $j$th subsample $\mathcal{S}_{j,k_n^*,r_n}$. Since these subsamples are separated by at least one block of size $r_n$ and since $\rho_{n,r_n}(u_{r_n}) = o(1/k_n) = o(1/q_n)$, a characteristic-function argument shows that the asymptotic distribution of $q_n^{-1/2}\sum_{j=1}^{q_n}\xi_{n,j}$ is the same as if the variables $\xi_{n,1},\ldots,\xi_{n,q_n}$ were independent.

We apply the Lindeberg-Feller central limit theorem with Lyapounov's condition. By Propositions 3.1, 3.2 and 3.3 applied for sample size $n^* = k^* r$ together with the fact that $q_n/k_n \sim 1/k_n^*$ and $q_n k_n/(n-r_n+1)^2 \sim 1/(r_n^2 k_n^*)$, we have

$$\operatorname{var}(\xi_{n,1}) \to \boldsymbol{a}^\top \Sigma \boldsymbol{a}, \qquad n \to \infty.$$

Lyapounov's condition finally requires that there exists $\delta > 0$ such that

$$\frac{q_n}{(\sqrt{q_n})^{2+\delta}} \operatorname{E}[|\xi_{n,1}|^{2+\delta}] \to 0, \qquad n \to \infty.$$

We will show that Lyapounov's condition holds for $2+\delta = p$ with $p$ as in (4.1).

For $\mu \geqslant 1$, integer $m \geqslant 1$ and real numbers $c_1,\ldots,c_m$, we have $|\sum_{i=1}^m c_i|^\mu \leqslant m^{\mu-1}\sum_{i=1}^m |c_i|^\mu$ (proof by Jensen's inequality). As a consequence, it is sufficient to show that for $v \in \{1,2,3\}$,

$$q_n^{-\delta/2} \operatorname{E}[|\xi_{n,1,v}|^{2+\delta}] \to 0, \qquad n \to \infty.$$

Above we have chosen $k_n^*$ in such a way that $k_n^* \to \infty$ and $k_n^* = o(k_n)$ as $n \to \infty$. Now we reinforce the latter requirement to

$$k_n^* = o\left(k_n^{\frac{\delta}{2(1+\delta)}}\right), \qquad n \to \infty. \tag{B.5}$$

Then for $v = 2$, since $|\bar{I}_{t,r}^{\mathrm{sl}}| \leqslant 1$, we simply have

$$\begin{aligned} q_n^{-\delta/2}\,\mathrm{E}[|\xi_{n,1,v}|^{2+\delta}] &= O\left(q_n^{-\delta/2}\left(\frac{\sqrt{q_n k_n}}{n}\right)^{2+\delta}(k_n^* r)^{2+\delta}\right) \\ &= O\left(\frac{q_n (k_n^*)^{2+\delta}}{k_n^{1+\delta/2}}\right) = O\left(\frac{(k_n^*)^{1+\delta}}{k_n^{\delta/2}}\right) = o(1). \end{aligned}$$

For $v \in \{1,3\}$, we have to proceed a little differently. Let $\zeta_{i,r}$ be equal to $\bar{I}_{i,r}^{\mathrm{dj}}$ if $v = 1$ and $\bar{N}_{i,r}$ if $v = 3$. Then

$$\begin{aligned} q_n^{-\delta/2}\,\mathrm{E}[|\xi_{n,1,v}|^{2+\delta}] &= q_n^{-\delta/2}(q_n/k_n)^{1+\delta/2}\,\mathrm{E}\left[\left|\textstyle\sum_{i=1}^{k_n^*}\zeta_{i,r_n}\right|^{2+\delta}\right] \\ &\leqslant \frac{q_n}{k_n^{1+\delta/2}}(k_n^*)^{1+\delta}\,\mathrm{E}\left[\textstyle\sum_{i=1}^{k_n^*}|\zeta_{i,r_n}|^{2+\delta}\right] \\ &= \frac{q_n (k_n^*)^{2+\delta}}{k_n^{1+\delta/2}}\,\mathrm{E}[|\zeta_{1,r_n}|^{2+\delta}]. \end{aligned}$$

Now $\mathrm{E}[|\zeta_{1,r_n}|^{2+\delta}] = O(1)$; for $v = 1$ this is obvious and for $v = 3$ this follows by condition (4.1). Again, requirement (B.5) on $k_n^*$ ensures that the right-hand side of the previous display is $o(1)$ as $n \to \infty$. □

## Appendix C: Proof of Proposition 5.1

In view of Theorem 4.1 and the continuous mapping theorem, we only need to show that $\hat{\sigma}^2_{n,r_n}(u_n)$ is weakly consistent for $\mathrm{var}(N)$. It is not difficult to see that, as $n \to \infty$,

$$\bar{N}_{n,r_n}(u_{r_n}) = \hat{\tau}_{n,r_n}(u_{r_n}) + o_p(1) = \tau_{r_n}(u_{r_n}) + o_p(1).$$

Define

$$\bar{\sigma}^2_{n,r}(u) := \frac{1}{n-r+1}\sum_{i=0}^{n-r}\big(N_{i,i+r}(u) - \tau_r(u)\big)^2.$$

Then as $n \to \infty$,

$$\bar{\sigma}^2_{n,r_n}(u_{r_n}) - \frac{n-2r_n+1}{n-r_n+1}\hat{\sigma}^2_{n,r_n}(u_{r_n}) = \big(\bar{N}_{n,r_n}(u_{r_n}) - \tau_{r_n}(u_{r_n})\big)^2 = o_p(1).$$

Hence we can proceed with $\bar{\sigma}^2_{n,r}(u)$ rather than with $\hat{\sigma}^2_{n,r}(u)$. By stationarity and the fact that $\mathrm{E}[N_r(u)] = r\bar{F}(u) = \tau_r(u_r)$, we have

$$\mathrm{E}[\bar{\sigma}^2_{n,r}(u)] = \mathrm{var}\big(N_r(u)\big) = \sigma_r^2(u).$$

By uniform integrability, the moment condition (4.1) implies $\sigma_r^2(u_r) \to \alpha m_2$ as $r \to \infty$ (1, Corollary p.338). Hence we only need to show that $\operatorname{var}\big(\bar\sigma^2_{n,r_n}(u_{r_n})\big) \to 0$ as $n \to \infty$.

Put $\xi_{n,i} := \{N_{i,i+r_n}(u_{r_n}) - \tau_{r_n}(u_{r_n})\}^2$. Then

$$\bar\sigma^2_{n,r_n}(u_{r_n}) = (n - r_n + 1)^{-1} \sum_{i=0}^{n-r_n} \xi_{n,i},$$

and thus, writing $\kappa_n := \operatorname{var}(\xi_{n,i})$,

$$\begin{aligned}
&\operatorname{var}\big(\bar\sigma^2_{n,r_n}(u_{r_n})\big) \\
&= \frac{1}{(n-r_n+1)^2} \sum_{i=0}^{n-r_n} \operatorname{var}(\xi_{n,i}) + \frac{2}{(n-r_n+1)^2} \sum\sum_{0 \leqslant i < j \leqslant n-r_n} \operatorname{cov}(\xi_{n,i}, \xi_{n,j}) \\
&= \frac{\kappa_n}{n-r_n+1} + \frac{2}{(n-r_n+1)^2} \sum_{s=1}^{n-r_n} (n - r_n + 1 - s) \operatorname{cov}(\xi_{n,0}, \xi_{n,s}).
\end{aligned}$$

The first term on the right-hand side converges to zero as $n \to \infty$ since $\kappa_n = O(1)$ by (5.2). The second term can be written as a sum $I_n + II_n$ according to whether $s \leqslant r_n$ or not. The first term can be bounded by

$$|I_n| \leqslant \frac{2r_n}{n-r_n+1}\kappa_n \to 0, \qquad n \to \infty.$$

The second term is

$$|II_n| \leqslant \frac{2}{n-r_n+1}\kappa_n \sum_{l=1}^{n} \rho_{n,l}(u_{r_n}),$$

which, in view of (3.1), converges as $n \to \infty$ to zero as well. □